\documentclass[12pt]{article}
\usepackage{amsmath,fullpage,amsthm,amssymb,xcolor,graphicx}

\newtheorem{thm}{Theorem}[section]
\newtheorem{corollary}{Corollary}[section]

\newtheorem{definition}{Definition}[section]
\newtheorem{remark}{Remark}[section]
\newtheorem{example}{Example}[section]

\title{On spectrally optimal duals of frames generated by graphs}
\author{Deepshikha \thanks{Department of Mathematics, Shyampur Siddheswari Mahavidyalaya, University of Calcutta, West Bengal 711312, India. Email: dpmmehra@gmail.com }\  \and Aniruddha Samanta \thanks{Theoretical Statistics and Mathematics Unit, Indian Statistical Institute, Kolkata-700108, India. Email: aniruddha.sam@gmail.com} 
}
\date{\today}
\begin{document}
	\maketitle
	\baselineskip=0.25in
	\begin{abstract}
Recently, the concept of frames generated by graphs has been introduced in \cite{D}. In this paper, we study spectrally optimal dual frames of frames generated by graphs. We show that if the frame is generated by a connected graph, then its canonical dual frame is the unique spectrally optimal dual frame for $1$-erasure and $2$-erasures. Further, we show that the canonical dual frames of frames generated by disconnected graphs are non-unique spectrally optimal dual frames for $1$-erasure and $2$-erasures. 
	\end{abstract}

	
	
	\section{Introduction}

 In order to solve certain complex problems in the non-harmonic Fourier series, Duffin and Schaeffer \cite{DS} developed the idea of frames. Frames are redundant sequences of vectors that span a vector space. The redundancy of frames makes them very useful in many applications, including data and image processing, signal transmission, image segmentation, coding theory, etc. Due to redundancy, the data transmitted using frames can be reconstructed with minimal error if erasures occur during the transmission. \\[1mm] 

 Finite frames are widely used in frame theory due to their significant relevance in applications. In \cite{D}, authors introduced a new class of finite frames using simple graphs and called them frames generated by graphs. They established some graph properties using frame properties and vice-versa. Consequently, it provides a beautiful connection between the three classical branches of mathematics, namely, frame theory, graph theory, and matrix theory. Further, frames generated by graphs have some special properties, and some of them are stated in the preliminaries section of this article.\\[-4mm]

 The problem of erasures commonly occurs due to transmission losses or disturbances during transmission, see \cite{BP,CK}. The problem of erasures can be handled in two ways. One way is to find a frame that can minimize the maximal error in reconstruction \cite{HP}. Another way is to find an optimal dual frame for a pre-selected frame in order to minimize the maximal error of reconstruction, see \cite{LH1}. In this paper, we focus on the second approach. A dual frame of a pre-selected frame that minimizes the maximal error of reconstruction is called an optimal dual frame. Various measurements are used to measure the error operator. If operator norm is used, then we use the term optimal dual frame. If the spectral radius is used to measure the error operator, then dual frames that minimize the maximum error are called spectrally optimal dual frames (or, in short, $SOD$-frames). In \cite{DS2}, authors studied optimal dual frames of frames generated by graphs. In this paper, we study the spectrally optimal dual frame of frames generated by graphs. We establish some interesting properties of $SOD$-frames for $1$-erasure and for $2$-erasures. \\[-4mm] 
		
	The outline of this paper is as follows. In Section 2, we provide the basic definitions and results of frame theory, spectrally optimal dual frames, spectral graph theory, and frames generated by graphs required in the rest of the article. In Section 3, we study spectrally optimal dual frames of frames generated by graphs for $1$-erasure. We show that the canonical dual frames of frames generated by connected graphs are the unique spectrally optimal dual frames for $1$-erasure. Further, if a frame is generated by a disconnected graph, then its canonical dual is a non-unique $SOD$-frame for $1$-erasure. Section 3 deals with spectrally optimal dual frames for $2$-erasures. First, we compute the spectral radius of the error operator for $2$-erasures of frames generated by graphs. We establish that canonical duals of frames generated by graphs are $SOD$-frames for $2$-erasures. Further, we show that they are non-unique $SOD$-frames for $2$-erasures if the frame is generated by a disconnected graph.
	
	\section{Preliminaries}
Throughout the rest of the article, we consider finite frames and $\mathbb{C}^k$ is a $k$-dimensional Hilbert space. We use the notation $[n]$ to denote the set $\{1,2,\ldots,n\}$ and $M=[m_{i,j}]_{i\in [p], j\in [q]}$ to denote the matrix $M$ of order $p\times q$ where $m_{i,j}$ is the $(i,j)$-th entry of $M$. Also, $diag(\alpha_1,\alpha_2,\ldots,\alpha_n)$ denotes the diagonal matrix of order $n$ with respective diagonal entries $\alpha_1,\alpha_2,\ldots,\alpha_n$. 

  A sequence of vectors $\Phi:=\{\phi_i\}_{i\in[n]}$ in $\mathbb{C}^k$ is called a \emph{frame} for $\mathbb{C}^k$ if there exist constants $0<A_o\leq B_o$ such that for any $f\in\mathbb{C}^k$, we have
\begin{align*}
A_o\|f\|^2\leq\sum_{i\in[n]}|\langle f, \phi_i\rangle|^2\leq B_o\|f\|^2.
\end{align*}
The constants $A_o$ and $B_o$ are called the \emph{lower} and \emph{upper frame bounds} of $\Phi$, respectively. For a frame $\Phi=\{\phi_i\}_{i\in[n]}$, the \emph{analysis operator} $T_{\Phi}:\mathbb{C}^k\rightarrow\mathbb{C}^n$ is defined as $T_{\Phi}(f)=\{\langle f,\phi_i\rangle\}_{i\in[n]}$. The analysis operator $T_{\Phi}$ is linear and bounded. The adjoint of $T_{\Phi}$, $T^*_{\Phi}:\mathbb{C}^n\rightarrow\mathbb{C}^k$, is called a \emph{synthesis operator} and defined by $T_{\Phi}(\{\alpha_i\}_{i\in[n]})=\sum\limits_{i\in[n]}\alpha_i \phi_i$. The canonical matrix representation of the synthesis operator is a $k\times n$ matrix such that its $i^{\text{th}}$ column is the vector $\phi_i$ that is  $[T^*_{\Phi}]=[\phi_1\,\,\phi_2\,\,\cdots\,\,\phi_n]$. The \emph{frame operator} $S_{\Phi}:\mathbb{C}^k\rightarrow\mathbb{C}^k$ is the composition of the synthesis operator and analysis operator such that $S_{\Phi}(f)=\sum\limits_{i\in[n]}\langle f,\phi_i\rangle\phi_i$. The operator $S_{\Phi}$ is bounded, linear, self-adjoint, positive, and invertible. The composition $T_{\Phi}T_{\Phi}^*:\mathbb{C}^n\rightarrow\mathbb{C}^n$ is called the \emph{Gramian operator} of $\Phi$. The canonical matrix representation of the Gramian operator is said to be the \emph{Gramian matrix}, defined by $\mathcal{G}_{\Phi}=[\langle \phi_j,\phi_i\rangle]_{i,j\in[n]}$ where $\langle \phi_j,\phi_i\rangle$ is the $(i,j)$-th entry of $\mathcal{G}_{\Phi}$. A frame $\Psi=\{\psi_i\}_{i\in[n]}$ of $\mathbb{C}^k$ is called a \emph{dual frame} of $\Phi$ if for any $f\in\mathbb{C}^k$, $f=\sum\limits_{i\in[n]}\langle f,\phi_i\rangle\psi_i=\sum\limits_{i\in[n]}\langle f,\psi_i\rangle\phi_i$. The sequence $\{S_{\Phi}^{-1}\phi_i\}_{i\in[n]}$ is a also a dual frame of $\Phi$ and it is called the \emph{canonical dual frame} of $\Phi$. We denote the set of all dual frames of a frame $\Phi$ by $D_{\Phi}$. The dual frames which are not canonical dual are called \emph{alternate dual frames}. If $\Psi=\{\psi_i\}_{i\in[n]}\in D_{\Phi}$ then there exist a sequence $\{\mu_i\}_{i\in[n]}\subset\mathbb{C}^k$ such that $\psi_i=S_{\Phi}^{-1}\phi_i+\mu_i$ for $i\in[n]$ and $\sum\limits_{i\in[n]}\langle f,\mu_i\rangle \phi_i=0$ for all $f\in\mathbb{C}^k$. We say two frames $\Phi_1=\{\phi_{1,i}\}_{i\in[n]}$ and $\Phi_2=\{\phi_{2,i}\}_{i\in[n]}$ are \emph{unitary equivalent} if there exists a unitary operator $V:\mathbb{C}^k\rightarrow\mathbb{C}^k$ such that $V(\phi_{1,i})=\psi_{1,i}$ for all $i\in[n]$. Refer to \cite{BL, CK1, C2, DGM} for the general theory of frames and its applications.

Suppose $\Phi=\{\phi_i\}_{i\in[n]}$ is a frame for $\mathbb{C}^k$ and $\Psi=\{\psi_i\}_{i\in[n]}$ is a dual frame of $\Phi$. Then for $\Lambda\subset[n]$, the error operator $E_{\Phi,\Psi,\Lambda}:\mathbb{C}^k\rightarrow\mathbb{C}^k$ is defined by 
	\begin{align*}
		E_{\Phi,\Psi,\Lambda}(f)=\sum_{i\in\Lambda}\langle f,\phi_i\rangle \psi_i.
	\end{align*}
For $1\leq r<n$, let $\rho^{(r)}_{\Phi,\Psi}:=\max\{\rho(E_{\Phi,\Psi,\Lambda}):\Lambda\subset[n] \text{ and }|\Lambda|=r\}$, where $\rho(E_{\Phi,\Psi,\Lambda})$ is the spectral radius of $E_{\Phi,\Psi,\Lambda}$. A dual frame $\Psi$ of the frame $\Phi$ is called a \emph{spectrally optimal dual frame} (or simply \emph{$SOD$-frame})  for $1$-erasure if $\rho^{(1)}_{\Phi}:=\min\left\{\rho^{(1)}_{\Phi,\widetilde{\Psi}}:\widetilde{\Psi}\in D_{\Phi}\right\}=\rho^{(1)}_{\Phi,\Psi}$. Spectrally optimal dual frames for $r$-erasures are defined inductively. For $1<r<n$, a dual frame $\Psi=\{\psi_i\}_{i\in[n]}$ of a frame $\Phi=\{\phi_i\}_{i\in[n]}$ is called \emph{spectrally optimal dual frame for $r$-erasures} if $\Psi$ is an spectrally optimal dual frame for $(r-1)$-erasures and $\rho^{(r)}_{\Phi}:=\min\{\rho^{(r)}_{\Phi,\widetilde{\Psi}}:\widetilde{\Psi}\in D_{\Phi}\}=\rho^{(r)}_{\Phi,\Psi}$. We denote the set of all $SOD$-frame of $\Phi$ for $r$-erasures by $SOD_{\Phi}^r$.
	  
\begin{remark}\cite{PHM}
    If $\Phi=\{\phi_i\}_{i\in[n]}$ is a frame with dual frame $\Psi=\{\psi_i\}_{i\in[n]}$, then $\rho^{(1)}_{\Phi,\Psi}=\max\{|\langle \psi_i,\phi_i\rangle|:i\in[n]\}$.
\end{remark} 
\begin{thm}\cite{PHM}\label{thmPHM}
		Suppose $\Phi=\{\phi_i\}_{i\in[n]}$ is a frame for $\mathbb{C}^k$. Then the following are equivalent: 
		\begin{enumerate}
			\item There exists a dual frame $\Psi=\{\psi_i\}_{i\in[n]}$ of $\Phi$ such that $\rho^{(1)}_{\Phi,\Psi}=\rho^{(1)}_\Phi=\frac{k}{n}$.
			\item There exists a dual frame $\Psi=\{\psi_i\}_{i\in[n]}$ of $\Phi$ such that $|\langle \psi_i,\phi_i\rangle|=\frac{k}{n}$ for all $i\in[n]$.
		\end{enumerate}
		\end{thm}
	
	\begin{thm}\label{thm2}\cite{NAS}
		Suppose $\Phi=\{\phi_i\}_{i\in[n]}$ is a frame for $\mathbb{C}^k$ and $V:\mathbb{C}^k\rightarrow\mathbb{C}^k$ is an invertible operator. Then a dual frame $\Psi=\{\psi_i\}_{i\in[n]}$ is (unique) spectrally optimal dual frame of $\Phi$ for $r$-erasures if and only if $V\Psi$ is (unique) spectrally optimal dual frame of $(V^{-1})^*\Phi$ for $r$-erasures.
	\end{thm}
	For the detailed study of optimal dual frames and spectrally optimal dual frames, readers may refer to \cite{DS1, LH, LH1, PHM}.

Graph theory is one of the classical branches of mathematics. Along with the deep theoretical aspect of graph theory, it has many real-life applications, which include measuring the connectivity of communication networks, hierarchical clustering, ranking of hyperlinks in search engines, image segmentation, and so on, see \cite{NA, SGT3, SGT2, SGT1}. Let $ G=(V(G),E(G)) $ be a simple graph with the vertex set $ V(G)=\{ v_1, v_2, \dots, v_n\}$  and edge set $ E(G) $. The vertices $ v_i $ and $ v_j $ are said to be \emph{adjacent} if they are connected by an edge and are denoted by  $ v_i \sim v_j $. The number of vertices adjacent to $ v_i $ is known as the \emph{degree of the vertex} $ v_i $ and is denoted by $ d(v_i)$ or simply by $ d_i $. The \emph{degree matrix} of a graph $G$ on $n$ vertices, denoted by $D(G):=diag(d_1,d_2,\ldots,d_n)$. The \emph{adjacency matrix} of $ G $ on $ n $ vertices is $ A(G):=(a_{i,j})_{n\times n} $, where
	\begin{align*}
	a_{i,j}=\begin{cases}
		1, \ \text{for }  v_i\sim v_j \text{ and } i\neq j \\
		0, \text{ elsewhere.}
	\end{cases}
\end{align*}
The matrix $ A(G) $ is simply written as $ A $, if the graph $G$ is clear from the context. The \emph{Laplacian matrix} of a graph $G$ is defined as $L(G):=D(G)-A(G)$. We simply write $ L(G) $ as $ L $ if $ G $ is clear from the context.  The Laplacian matrix is a positive semi-definite matrix. For a simple graph $G$ on $n$ vertices with $k$ components, the rank of $L(G)$ is $n-k$. For a detailed study of matrices associated with simple graphs and other classes of graphs, readers may see \cite{Bapat, SGT6, SGT7, SGT8, SGT9, SGT10} and references therein. The transpose of a matrix $ M $ is denoted by $ M^t $. If $ M_1 $ and $ M_2 $ are matrices, then $ M_1 \oplus M_2 $ denotes the block matrix $\left[\begin{array}{cc}
		M_1 & \boldsymbol{0}  \\
		\boldsymbol{0} & M_2
	\end{array}\right]$.

In \cite{D}, authors used Laplacian matrices of graphs to generate finite frames. First, let us see the definition of $L_G(n,k)$-frames.
	\begin{definition}
		Let $G$ be a simple graph on  $n$ vertices with $n-k$ components. Suppose $L$ is the Laplacian matrix of $G$ such that $L=M\, diag(\lambda_1,\lambda_2,\ldots,\lambda_{k},0,\ldots,0)\,M^*$. If $\{e_i\}_{i\in[n]}$ is the standard canonical orthonormal basis of $\mathbb{C}^n$ and $B=diag\left(\sqrt{\lambda_1},\ldots,\sqrt{\lambda_k}\right)M_1^*$, where $M_1$ is a submatrix of $M$ formed by the first $k$ columns then $\{B(e_i)\}_{i\in[n]}$ is a frame for $\mathbb{C}^k$ and it is called an $L_G(n,k)$-frame for $\mathbb{C}^k$.
	\end{definition}
	Next, we have the definition of frames generated by graphs. 
	\begin{definition}
		Let $\Phi=\{\phi_i\}_{i\in [n]}$ be a frame for $\mathbb{C}^k$ with Gramian matrix $\mathcal{G}_{\Phi}$. If $G$ is a graph with Laplacian matrix $L$ such that $\mathcal{G}_{\Phi}=L$, then $\Phi$ is called a frame generated by the graph $G$ for $\mathbb{C}^k$. In short, we call $\Phi$ as a $G(n,k)$-frame for $\mathbb{C}^k$.
	\end{definition}
	In \cite{D}, the authors have shown that every $L_G(n,k)$-frame is a $G(n,k)$-frame. Let us see some results related to frames generated by graphs.
	
	\begin{thm}[\cite{D}]\label{thm2.3}
		If $G$ is a graph, then frames generated by $G$ are unitary equivalent.
	\end{thm}
	The next theorem shows that the frame operator of any $L_G(n,k)$-frame is a diagonal matrix. 
	\begin{thm}[\cite{D}]\label{thm2.7}
		If $G$ is a graph and $\Phi$ is an $L_G(n,k)$-frame then the frame operator $S_{\Phi}$ of the frame $\Phi$ is a diagonal matrix. Further, if $\lambda_1,\lambda_2,\ldots,\lambda_k$ are the non-zero eigenvalues of the Laplacian matrix of $G$ then $S_{\Phi}=diag(\lambda_1,\lambda_2,\ldots,\lambda_k)$.
	\end{thm}
	The following theorem gives the family of dual frames of the frames generated by graphs.
	\begin{thm}\label{thm2.5}
		Let $G$ be a simple graph with components $G_1,G_2,\ldots, G_{m}$ having vertex sets $\{v_1,v_2,\ldots,v_{n_1}\}, \{v_{n_1+1}, v_{n_1+2}, \ldots, v_{n_2}\}, \ldots,\{v_{n_{m-1}+1}, v_{n_{m-1}+2}, \ldots, v_{n_m}(=v_n)\}$, respectively. If $\Phi=\{\phi_i\}_{i\in[n]}$ is a $G(n,n-m)$-frame with the frame operator $S_{\Phi}$, then any dual frame of $\Phi$ is of the form $\{S_{\Phi}^{-1}\phi_i+\nu_1\}_{i=1}^{n_1}\bigcup\{S_{\Phi}^{-1}\phi_i+\nu_2\}_{i=n_1+1}^{n_2}\bigcup\cdots\bigcup\{S_{\Phi}^{-1}\phi_i+\nu_m\}_{i=n_{m-1}+1}^{n}$ where $\nu_1,\nu_2,\ldots,\nu_m$ are arbitrary vectors in $\mathbb{C}^{n-m}$.
	\end{thm}
 
	\section{On spectrally optimal dual frame for $1$-erasure}
  We begin this article by showing that if the canonical dual frame of a frame generated by a graph $G$ is a spectrally optimal dual frame for $r$-erasures, then the canonical dual of any frame generated by $G$ is a spectrally optimal dual frame for $r$-erasures.
\begin{thm}\label{thm4.2}
	Let $G$ be any graph with $n$ vertices. If $\Phi_1=\{\phi_{1,i}\}_{i\in[n]}$ and $\Phi_2=\{\phi_{2,i}\}_{i\in[n]}$ are $G(n,k)$-frames for $\mathbb{C}^k$, then the canonical dual of $\Phi_1$ is (unique) spectrally optimal dual frame of $\Phi_1$ for $r$-erasures if and only if the canonical dual of $\Phi_2$ is (unique) spectrally optimal dual frame of $\Phi_2$ for $r$-erasures. 
\end{thm}	
\proof
Since $\Phi_1$ and $\Phi_2$ are $G(n,k)$-frames, by Theorem \ref{thm2.3}, there is a unitary operator $V:\mathbb{C}^k\rightarrow\mathbb{C}^k$ such that $V(\phi_{1,i})=\phi_{2,i}$ for all $i\in[n]$. Let $S_{\Phi_1}$ and $S_{\Phi_2}$ are the frame operators of $\Phi_1$ and $\Phi_2$, respectively. Then, for any $f\in\mathbb{C}^k$, we have
\begin{align*}
S_{\Phi_2}(f)&=\sum\limits_{i\in[n]}\langle f,\phi_{2,i}\rangle \phi_{2,i}\\
&=\sum\limits_{i\in[n]}\langle f,V\phi_{1,i}\rangle V\phi_{1,i}\\
&=VS_{\Phi_1}V^*(f).
\end{align*}
Thus, $S_{\Phi_2}=VS_{\Phi_1}V^*$ that is $S_{\Phi_1}^{-1}=V^*S_{\Phi_2}^{-1}V$. Then, by Theorem \ref{thm2}, $S_{\Phi_1}^{-1}\Phi_1$ is (unique) $SOD$-frame of $\Phi_1$ for $r$-erasures if and only if $VS_{\Phi_1}^{-1}\Phi_1=VV^*S_{\Phi_2}^{-1}VV^*\Phi_2=S_{\Phi_2}^{-1}\Phi_2$ is (unique) $SOD$-frame of $(V^{-1})^*\Phi_1=V\Phi_1=\Phi_2$ for $r$-erasures.
\endproof
In the following theorem, we show that the canonical dual of a frame generated by the connected graph is the unique $SOD$-frame for $1$-erasure.
\begin{thm}\label{thm4.1}
Suppose $G$ is a connected graph with $n$ vertices. If $\Phi=\{\phi_i\}_{i\in[n]}$ is an $L_G(n,n-1)$-frame for $\mathbb{C}^{n-1}$, then canonical dual frame of $\Phi$ is the unique spectrally optimal dual frame of $\Phi$ for $1$-erasure. Further, $\rho^{(1)}_{\Phi}=\frac{n-1}{n}$.
\end{thm}	
\proof
Let $L$ be the Laplacian matrix of $G$ with eigenvalues $\lambda_1,\lambda_2,\ldots,\lambda_{n-1},0$. Then there exist an orthogonal matrix $M$ consisting of eigenvectors of $L$ and $D=diag(\lambda_1,\lambda_2,\ldots ,\lambda_{n-1},0)$ such that $L=MDM^*$ and  $\phi_i=D_1M_1^*(e_i)$ where $D_1=diag(\sqrt{\lambda_1},\sqrt{\lambda_2},\ldots ,\sqrt{\lambda_{n-1}})$, $M_1$ is the matrix obtained from $M$ by taking the first $n-1$ columns of $M$ and $\{e_i\}_{i\in[n]}$ is the standard canonical orthonormal basis of $\mathbb{C}^n$. Note that $[1\,1\,\cdots\,1]^t$ is an eigenvector of $L$ corresponding to the eigenvalue $0$. Hence, the last column of the matrix $M$ is $\left[\frac{1}{\sqrt{n}}\, \frac{1}{\sqrt{n}}\,\cdots\,\frac{1}{\sqrt{n}}\right]^t$. Suppose $S_{\Phi}$ is the frame operator of $\Phi$. Then, by Theorem \ref{thm2.7}, $S_{\Phi}=diag(\lambda_1,\lambda_2,\ldots,\lambda_{n-1})$ and hence $S_{\Phi}^{1/2}=diag(\sqrt{\lambda_1},\sqrt{\lambda_2},\ldots,\sqrt{\lambda_{n-1}})=D_1$. Thus $S_{\Phi}^{-1/2}=D_1^{-1}$. Then we have
\begin{align*}
|\langle S_{\Phi}^{-1}\phi_i,\phi_i\rangle|=\|S_{\Phi}^{-1/2}\phi_i\|^2=\|D_1^{-1}D_1M_1^*(e_i)\|^2=\|M_1^*(e_i)\|^2=1-\frac{1}{n}.
\end{align*}
Thus, $|\langle S_{\Phi}^{-1}\phi_i,\phi_i\rangle|=\frac{n-1}{n}$ for all $i\in[n]$. Hence, by Theorem \ref{thmPHM}, $\rho_{\Phi,S_{\Phi}^{-1}\Phi}^{(1)}=\rho_{\Phi}^{(1)}=\frac{n-1}{n}$. Therefore, the canonical dual $S_{\Phi}^{-1}\Phi$ is an $SOD$-frame of $\Phi$ for $1$-erasure. Let $\Psi=\{\psi_i\}_{i\in[n]}$ be a spectrally optimal dual frame of $\Phi$ for $1$-erasure. Then, by Theorem \ref{thm2.5}, there exists a vector $\mu\in\mathbb{C}^{n-1}$ such that $\psi_i=S_{\Phi}^{-1}\phi_i+\mu$ for all $i\in[n]$. Then $\sum\limits_{i\in[n]}\langle f,\mu\rangle \phi_i=0$ for all $f\in\mathbb{C}^{n-1}$. Thus, $\sum\limits_{i\in[n]}\langle \mu,\phi_i\rangle = 0$ and we have
\begin{align*}
	\rho_{\Phi,\Psi}^{(1)}&=\max\{|\langle S_{\Phi}^{-1}\phi_i+\mu,\phi_i\rangle|:i\in[n]\}\\
	&=\max\{|\|S_{\Phi}^{-1/2}\phi_i\|^2+\langle\mu,\phi_i\rangle|:i\in[n]\}\\
			&=\max\left\{\left|\frac{n-1}{n}+\langle\mu,\phi_i\rangle\right|:i\in[n]\right\}\\
			&=\max\left\{\sqrt{\left|\frac{n-1}{n}+Re\langle\mu,\phi_i\rangle\right|^2+|Im\langle\mu,\phi_i\rangle|^2}:i\in[n]\right\}\\
			&=\frac{n-1}{n}.
\end{align*}
Thus, $\left|\frac{n-1}{n}+Re\langle\mu,\phi_i\rangle\right|\leq\frac{n-1}{n}$ for all $i\in[n]$. Hence, $Re\langle\mu,\phi_i\rangle\leq 0$ for all $i\in[n]$. Since $\sum\limits_{i\in[n]}\langle \mu,\phi_i\rangle = 0$ , $Re\langle \mu,\phi_i\rangle=0$ for $i\in[n]$. Thus, $\rho_{\Phi,\Psi}^{(1)}=\max\left\{\sqrt{\left|\frac{n-1}{n}\right|^2+|Im\langle\mu,\phi_i\rangle|^2}:i\in[n]\right\}\break =\frac{n-1}{n}$. This gives that $|Im\langle\mu,\phi_i\rangle|^2=0$ for all $i\in[n]$. Hence, $Im\langle\mu,\phi_i\rangle=0$ for all $i\in[n]$. Then $\langle\mu,\phi_i\rangle=0$ for all $i\in[n]$ and hence $\mu=0$. Thus, $\Psi=S_{\Phi}^{-1}\Phi$ and $S_{\Phi}^{-1}\Phi$ is the unique spectrally optimal dual frame of $\Phi$ for $1$-erasure.
\endproof

\begin{corollary}
Suppose $G$ is a connected graph with $n$ vertices. If $\Phi=\{\phi_i\}_{i\in[n]}$ is a $G(n,n-1)$-frame for $\mathbb{C}^{n-1}$, then canonical dual frame of $\Phi$ is the unique spectrally optimal dual frame of $\Phi$ for $1$-erasure.	
\end{corollary}
\proof
Proof follows from Theorem \ref{thm4.2} and Theorem \ref{thm4.1}.
\endproof
In Theorem \ref{thm4.1}, we consider $G$ to be a connected graph. In the following theorem, we show that Theorem \ref{thm4.1} is partially true for disconnected graphs.
\begin{thm}\label{thm4.3}
	Suppose $G$ is a graph with $n$ vertices and $k>1$ connected components. If $\Phi=\{\phi_i\}_{i\in[n]}$ is an $L_G(n,n-k)$-frame for $\mathbb{C}^{n-k}$, then canonical dual frame of $\Phi$ is a non-unique spectrally optimal dual frame of $\Phi$ for $1$-erasure. 
\end{thm}	
\proof
Let $G_1,G_2,\ldots,G_k$ be the connected components of $G$ with $|V(G_i)|=n_i$ for $i\in[k]$. Let the vertex set of $G_i$ be $\{v_{l_{i-1}+1},v_{l_{i-1}+2},\ldots,v_{l_i}\}$ for $i\in[k]$ where $l_0=0$ and for $1\leq i\leq k$, $l_i=n_1+n_2+\cdots+n_i$. If $L_i$ is the Laplacian matrix of $G_i$, then the Laplacian matrix of $G$ is $L=L_1\oplus L_2\oplus\cdots\oplus L_k$. Suppose the the spectral decomposition of $L_i$ is $M_iD_iM_i^*$w that is $L_i=M_iD_iM_i^*$ where $D_i=diag(\lambda_1^i,\ldots,\lambda_{n_i-1}^i,0)$ and $M_i$ is an orthogonal matrix consist of eigenvectors of $L_i$. Then $L=MDM^*$ where $M=M_1\oplus\cdots\oplus M_k$ and $D=D_1\oplus\cdots\oplus D_k$. Let $\widetilde{D_i}=diag\left(\sqrt{\lambda_1^i},\ldots,\sqrt{\lambda_{n_i-1}^i}\right)$ and $\widetilde{M_i}$ be obtained from $M_i$ by taking the first $n_i-1$ columns for $i\in[k]$. Then $\widetilde{D_i}\widetilde{M_i}^*$ is the synthesis operator of frame say $\Phi_i=\{\phi_{i,j}\}_{j\in[n_i]}$ for $\mathbb{C}^{n_i-1}$. Let $S_{\Phi_i}$ be the frame operator of $\Phi_i$ for $i\in[k]$. 

Let $\widetilde{D}=\widetilde{D_1}\oplus\cdots\oplus\widetilde{D_k}$ and $\widetilde{M}
=\widetilde{M_1}\oplus\cdots\oplus\widetilde{M_k}$. Then $\widetilde{\Phi}=\{\widetilde{\phi}_i\}_{i\in[n]}=\{\widetilde{D}\widetilde{M}^*(e_i)\}_{i\in[n]}$ is an $L_G(n,n-k)$-frame for $\mathbb{C}^{n-k}$. Also, the synthesis operator $T_{\widetilde{\Phi}}^*$ of $\widetilde{\Phi}$ is
\begin{align*}
	[T_{\widetilde{\Phi}}^*]&=[T_{\Phi_1}^*]\oplus[T_{\Phi_2}^*]\oplus\cdots\oplus[T_{\Phi_k}^*]\\
	&=\left[\begin{array}{cccccccccc}
		\phi_{1,1}  &\cdots & \phi_{1,n_1} & 0 & \cdots & 0 & \cdots & 0 & \cdots & 0\\
		0  &\cdots & 0 & \phi_{2,1} & \cdots & \phi_{2,n_2} & \cdots & 0 & \cdots & 0\\
		\vdots  &\ddots & \vdots & \vdots & \ddots & \vdots & \ddots & \vdots & \ddots & \vdots \\
		0  &\cdots & 0 & 0 & \cdots & 0 & \cdots & \phi_{k,1} & \cdots & \phi_{k,n_k}
	\end{array}\right]\\
	&=\{\widetilde{\phi}_i\}_{i\in[n]}
\end{align*}

By Theorem \ref{thm4.2}, it is enough to prove that the canonical dual of $\widetilde{\Phi}$ is a non-unique spectrally optimal dual frame of $\widetilde{\Phi}$ for $1$-erasure. 

Let $S_{\widetilde{\Phi}}$ be the frame operator of $\widetilde{\Phi}$. Then $S_{\widetilde{\Phi}}=S_{\Phi_1}\oplus\cdots\oplus S_{\Phi_k}$. Suppose $\max\{|\langle S_{\widetilde{\Phi}}^{-1}\widetilde{\phi}_i,\widetilde{\phi}_i\rangle|:i\in[n]\}=|\langle S_{\widetilde{\Phi}}^{-1}\widetilde{\phi}_j,\widetilde{\phi}_j\rangle|$ where $j\in[l_{r-1}+1,\ldots,l_r]$ for some $r\in[k]$. Then $\rho^{(1)}_{\widetilde{\Phi},S_{\widetilde{\Phi}}^{-1}\widetilde{\Phi}}=|\langle S_{\widetilde{\Phi}}^{-1}\widetilde{\phi}_j,\widetilde{\phi}_j\rangle|=1-\frac{1}{n_r}=\frac{n_r-1}{n_r}$. By Theorem \ref{thm4.1}, $\rho^{(1)}_{\Phi_r}=\rho^{(1)}_{\Phi_r,S_{\Phi_r}^{-1}\Phi_r}=\frac{n_r-1}{n_r}$.

 Suppose $\Psi=\{\psi_i\}_{i\in[n]}$ is a dual frame of $\widetilde{\Phi}$. Since $\Psi$ is a dual frame, there exist vectors $\mu_1,\mu_2,\ldots,\mu_k$ in $\mathbb{C}^{n-k}$ such that $\psi_i=S_{\widetilde{\Phi}}^{-1}\widetilde{\phi}_i+\mu_j$  for $i\in\{l_{j-1}+1,l_{j-1}+2,\ldots,l_{j}\}$. Let $\mu_j=\left[\begin{array}{c}
      \nu_1 \\
      \nu_2\\
      \vdots\\
      \nu_k
 \end{array}\right]$ where $\nu_i\in\mathbb{C}^{n_i-1}$. Then by Theorem \ref{thm2.5}, $\Psi_r=\{S_{\Phi_r}^{-1}\phi_{r,i}+\nu_r\}_{i\in[n_r]}$ is a dual frame of $\Phi_r$. By Theorem \ref{thm4.1}, $\rho_{\Phi_r,\Psi_r}^{(1)}\geq\frac{n_r-1}{n_r}$. Then
 \begin{align*}
 	\rho^{(1)}_{\widetilde{\Phi},\Psi}&=\max\{|\langle \psi_i,\widetilde{\phi}_i\rangle|:i\in[n]\}\\
 	&\geq \max\{|\langle \psi_i,\widetilde{\phi}_i\rangle|:i\in\{l_{r-1}+1,\ldots,l_r\}\}\\
 	&=\max\{|\langle S_{\widetilde{\Phi}}^{-1}\widetilde{\phi}_i+\mu_j,\widetilde{\phi}_i\rangle|:i\in\{l_{r-1}+1,\ldots,l_r\}\}\\
 	&=\max\{|\langle S_{\Phi_r}^{-1}\widetilde{\phi}_{r,i}+\nu_r,\widetilde{\phi}_{r,i}\rangle|:i\in[n_r]\}\\
 	&=\rho^{(1)}_{\Phi_r,\Psi_r}\\
 	&\geq \frac{n_r-1}{n_r}\\
&=\rho^{(1)}_{\widetilde{\Phi},S_{\widetilde{\Phi}}^{-1}\widetilde{\Phi}}=\rho^{(1)}_{\widetilde{\Phi}}.
 \end{align*}
Therefore, $S_{\widetilde{\Phi}}^{-1}\widetilde{\Phi}$ is an $SOD$-frame of $\widetilde{\Phi}$ for $1$-erasure. Then by Theorem \ref{thm4.2}, $S_{\Phi}^{-1}\Phi$ is an $SOD$-frame of $\Phi$ for $1$-erasure.

Let $\mu=[\mu_i]_{i\in[n]}\in\mathbb{C}^{n-k}$ such that $\mu_i=\begin{cases}
    0, \text{ if } i\in[l_1]\\
    1, \text{ otherwise}\\
\end{cases}$. By Theorem \ref{thm2.5}, $\Psi_1=\{\psi_{1,i}\}_{i\in[n]}=\{S_{\widetilde{\Phi}}^{-1}\widetilde{\phi}_i+\mu\}_{i\in[l_1]}\cup\{S_{\widetilde{\Phi}}^{-1}\widetilde{\phi}_i\}_{i\in[n]\setminus[l_1]}$ is a dual frame of $\widetilde{\Phi}$. Note that $\Psi_1$ is not the canonical dual frame of $\widetilde{\Phi}$. We have
\begin{align*}
    \rho_{\widetilde{\Phi},\Psi_1}^{(1)}&=\max\{|\langle \psi_{1,i},\widetilde{\phi}_i\rangle|:i\in[n]\}\\
    &=\max(\{|\langle S_{\widetilde{\Phi}}^{-1}\widetilde{\phi}_i+\mu,\widetilde{\phi}_i\rangle|:i\in[l_1]\}\cup\{|\langle S_{\widetilde{\Phi}}^{-1}\widetilde{\phi}_i,\widetilde{\phi}_i\rangle|:i\in[n]\setminus[l_1]\})\\
    &=\max(\{|\langle S_{\widetilde{\Phi}}^{-1}\widetilde{\phi}_i,\widetilde{\phi}_i\rangle+\langle \mu,\widetilde{\phi}_i\rangle|:i\in[l_1]\}\cup\{|\langle S_{\widetilde{\Phi}}^{-1}\widetilde{\phi}_i,\widetilde{\phi}_i\rangle|:i\in[n]\setminus[l_1]\})\\
     &=\max(\{|\langle S_{\widetilde{\Phi}}^{-1}\widetilde{\phi}_i,\widetilde{\phi}_i\rangle:i\in[n]\}\cup\{|\langle S_{\widetilde{\Phi}}^{-1}\widetilde{\phi}_i,\widetilde{\phi}_i\rangle|:i\in[n]\})\\
     &=\rho_{\widetilde{\Phi},S_{\widetilde{\Phi}}^{-1}\widetilde{\Phi}}^{(1)}.
\end{align*}
Hence, $\Psi_1$ is an $SOD$-frame of $\widetilde{\Phi}$ for $1$-erasure. Therefore, by Theorem \ref{thm4.2}, the canonical dual of $\Phi$ is a non-unique $SOD$-frame of $\Phi$ for $1$-erasure.
\endproof

\begin{corollary}
	Suppose $G$ is a graph with $n$ vertices and $k>1$ connected components. If $\Phi=\{\phi_i\}_{i\in[n]}$ is a $G(n,n-k)$-frame for $\mathbb{C}^{n-k}$, then the canonical dual frame of $\Phi$ is a non-unique spectrally optimal dual frame of $\Phi$ for $1$-erasure.	
\end{corollary}

In the following example, we show explicitly that the canonical dual frames of frames generated by disconnected graphs are non-unique spectrally optimal dual frames for $1$-erasure.

\begin{example}\label{Exa1}

  {\em  Consider the graph $G$ given in Figure \ref{fig1}.
    \begin{figure}
			\begin{center}
				\includegraphics[scale= 0.65]{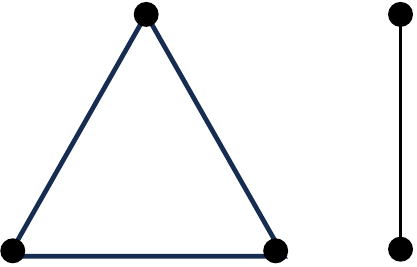}
				\caption{ Graph $ G $} \label{fig1}
			\end{center}
		\end{figure} The graph $G$ is a disconnected graph. The Laplacian matrix of $G$ is $L(G)=\left[\begin{array}{ccccc}
    2 & -1 & -1 & 0 & 0\\
    -1 & 2 & -1 & 0 & 0\\
    -1 & -1 & 2 & 0 & 0\\
    0 & 0 & 0 & 1 & -1\\
    0 & 0 & 0 & -1 & 1
    \end{array}\right]$.\\[3mm] 
    
    \noindent Consider the frame $\Phi=\{\phi_i\}_{i\in[5]}=\left\{\left[\begin{array}{c}
    \frac{\sqrt{3}}{\sqrt{2}} \\[2mm]
    \frac{-1}{\sqrt{2}}  \\[2mm]
    0 
    \end{array}\right], \left[\begin{array}{c}
    0 \\[2mm]
    \sqrt{2}  \\[2mm]
    0 
    \end{array}\right], \left[\begin{array}{c}
    \frac{-\sqrt{3}}{\sqrt{2}} \\[2mm]
    \frac{-1}{\sqrt{2}}  \\[2mm]
    0 
    \end{array}\right], \left[\begin{array}{c}
    0\\[2mm]
    0\\[2mm]
    1
    \end{array}\right], \left[\begin{array}{c}
    0\\[2mm]
    0\\[2mm]
    -1
    \end{array}\right]\right\}$ of $\mathbb{C}^3$.\\[3mm]
 Then the Gramian matrix of $\Phi$ is $\mathcal{G}_{\Phi}=\left[\begin{array}{ccccc}
    2 & -1 & -1 & 0 & 0\\
    -1 & 2 & -1 & 0 & 0\\
    -1 & -1 & 2 & 0 & 0\\
    0 & 0 & 0 & 1 & -1\\
    0 & 0 & 0 & -1 & 1
    \end{array}\right] = L(G)$. Thus, $\Phi$ is a \\[2mm] $G(5,3)$-frame for $\mathbb{C}^3$. Let $S_{\Phi}$ be the frame operator of $\Phi$. Then the canonical dual frame of \\[2mm] $\Phi$ is $S_{\Phi}^{-1}\Phi =\left\{\left[\begin{array}{c}
    \frac{1}{\sqrt{6}} \\[2mm]
    \frac{-1}{3\sqrt{2}}  \\[2mm]
    0 
    \end{array}\right], \left[\begin{array}{c}
    0 \\[2mm]
    \frac{2}{3\sqrt{2}}  \\[2mm]
    0 
    \end{array}\right], \left[\begin{array}{c}
    \frac{-1}{\sqrt{6}} \\[2mm]
    \frac{-1}{3\sqrt{2}}  \\[2mm]
    0 
    \end{array}\right], \left[\begin{array}{c}
    0\\[2mm]
    0\\[2mm]
    \frac{1}{2}
    \end{array}\right], \left[\begin{array}{c}
    0\\[2mm]
    0\\[2mm]
    \frac{-1}{2}
    \end{array}\right]\right\}$. By Theorem \ref{thm2.5}, the set of\\[2mm] dual frames of $\Phi$ is $$D_{\Phi}=\left\{\left\{\left[\begin{array}{c}
    \frac{1}{\sqrt{6}} + a\\[2mm]
    \frac{-1}{3\sqrt{2}} + b \\[2mm]
    c
    \end{array}\right], \left[\begin{array}{c}
    a \\[2mm]
    \frac{2}{3\sqrt{2}} + b\\[2mm]
    c 
    \end{array}\right], \left[\begin{array}{c}
    \frac{-1}{\sqrt{6}} + a \\[2mm]
    \frac{-1}{3\sqrt{2}} + b \\[2mm]
    c 
    \end{array}\right], \left[\begin{array}{c}
    e\\[2mm]
    f\\[2mm]
    \frac{1}{2} +g
    \end{array}\right], \left[\begin{array}{c}
    e\\[2mm]
    f\\[2mm]
    \frac{-1}{2} + g
    \end{array}\right]\right\}: a, b, c, e ,f , g \in\mathbb{C}\right\}$$.

    Consider $\rho^{(1)}_{\Phi, S_{\Phi}^{-1}\Phi}=\max\{|\langle S_{\Phi}^{-1}\phi_i,\phi_i\rangle|: i\in[5]\}=\max\left\{\frac{2}{3},\frac{1}{2}\right\}=\frac{2}{3}$. Let $\Psi=\{\psi_i\}_{i\in[5]}$ be any dual frame of $\Phi$. Then there exist scalars $a_o,b_o,c_o,d_o,e_o,f_o,g_o$ such that $\Psi=\{\psi_i\}_{i\in[5]}=\left\{\left[\begin{array}{c}
    \frac{1}{\sqrt{6}} + a_o\\[1mm]
    \frac{-1}{3\sqrt{2}} + b_o \\[1mm]
    c_o
    \end{array}\right], \left[\begin{array}{c}
    a_o \\[1mm]
    \frac{2}{3\sqrt{2}} + b_o\\[1mm]
    c_o 
    \end{array}\right], \left[\begin{array}{c}
    \frac{-1}{\sqrt{6}} + a_o \\[1mm]
    \frac{-1}{3\sqrt{2}} + b_o \\[1mm]
    c_o 
    \end{array}\right], \left[\begin{array}{c}
    e_o\\[1mm]
    f_o\\[1mm]
    \frac{1}{2} +g_o
    \end{array}\right], \left[\begin{array}{c}
    e_o\\[1mm]
    f_o\\[1mm]
    \frac{-1}{2} + g_o
    \end{array}\right]\right\}$. Then 
    \begin{align*}
    \rho^{(1)}_{\Phi, \Psi}&=\max\{|\langle \psi_i,\phi_i\rangle|: i\in[5]\}\\
    &=\max\left\{\frac{2}{3}+\frac{\sqrt{3}}{\sqrt{2}}a_o-\frac{1}{\sqrt{2}}b_o,\frac{2}{3}+\sqrt{2}b_o,\frac{2}{3}-\frac{\sqrt{3}}{\sqrt{2}}a_o-\frac{1}{\sqrt{2}}b_o,\frac{1}{2}+g_o,\frac{1}{2}-g_o\right\}\\
    &\geq\frac{2}{3}=\rho^{(1)}_{\Phi, S_{\Phi}^{-1}\Phi}.
    \end{align*}
    Thus, $S^{-1}\Phi$ is an $SOD$-frame of $\Phi$ for $1$-erasure. 
    Consider $$\Psi_1=\{\psi_{1,i}\}_{i\in[5]}=\left\{\left[\begin{array}{c}
    \frac{1}{\sqrt{6}}\\[1mm]
    \frac{-1}{3\sqrt{2}}\\[1mm]
    1
    \end{array}\right], \left[\begin{array}{c}
    0 \\[1mm]
    \frac{2}{3\sqrt{2}}\\[1mm]
    1
    \end{array}\right], \left[\begin{array}{c}
    \frac{-1}{\sqrt{6}} \\[1mm]
    \frac{-1}{3\sqrt{2}} \\[1mm]
    1
    \end{array}\right], \left[\begin{array}{c}
    0\\[1mm]
    0\\[1mm]
    \frac{1}{2} 
    \end{array}\right], \left[\begin{array}{c}
    0\\[1mm]
    0\\[1mm]
    \frac{-1}{2}
    \end{array}\right]\right\}.$$
    Then $\Psi_1$ is a dual frame of $\Phi$. Also, $\rho^{(1)}_{\Psi,\Phi}=\max\{|\langle \psi_{1,i},\phi_i\rangle|: i\in[5]\}=\frac{2}{3}$. Thus, $\Psi_1$ is an alternate dual frame of $\Phi$ and an $SOD$-frame of $\Phi$ for $1$-erasure. Thus, the canonical dual of $\Phi$ is a non-unique spectrally optimal dual of $\Phi$ for $1$-erasure.}
\end{example}

\section{On Spectrally optimal dual frames for $2$-erasures}
In this section, we study spectrally optimal dual frames of frames generated by graphs for $2$-erasures. In the previous section, we establish that frames generated by graphs are spectrally optimal dual frames for $1$-erasure. Here, we begin by measuring the error operator if canonical dual is used and $2$-erasures occur. 
\begin{thm}\label{thm5.1}
Suppose $G$ is a connected graph on $n$ vertices and $\Phi=\{\phi_i\}_{i\in[n]}$ is an $L_G(n,n-1)$-frame for $\mathbb{C}^{n-1}$. If $S_{\Phi}$ is the frame operator of $\Phi$ then $\rho_{\Phi,S_{\Phi}^{-1}\Phi}^{(2)}=1$.
\end{thm}
\proof
Let $L$ be the Laplacian matrix of $G$. Since $\Phi$ is an $L_G(n,n-1)$-frame, there exist an orthogonal matrix $M$ and diagonal matrix $D=diag(\lambda_1,\lambda_2,\ldots,\lambda_{n-1},0)$ such that $L=MDM^*$ and $\Phi=D_1M_1^*$ where $M_1$ is formed by taking the first $n-1$ columns of $M$ and $D_1=diag(\sqrt{\lambda_1},\sqrt{\lambda_2},\ldots,\sqrt{\lambda_{n-1}})$. Note that the last column of $M$ is $\left[\begin{array}{c}
     \frac{1}{\sqrt{n}}\\[1mm]
     \frac{1}{\sqrt{n}}\\[1mm]
     \vdots\\
     \frac{1}{\sqrt{n}}\\[1mm]
\end{array}\right]$. Then, by Theorem \ref{thm2.7}, $S_{\Phi}=diag(\lambda_1,\ldots,\lambda_{n-1})$. Let $\Lambda=\{i,j\}\subset[n]$ be arbitrary. Then, the error operator $E_{\Phi, S_{\Phi}^{-1}\Phi,\Lambda}:\mathbb{C}^k\rightarrow\mathbb{C}^k$ is defined by
$$E_{\Phi,S_{\Phi}^{-1}\Phi,\Lambda}(f)=\langle f, \phi_i\rangle S_{\Phi}^{-1}\phi_i+\langle f,\phi_j\rangle S_{\Phi}^{-1}\phi_j.$$ Consider 
\begin{align*}
    E_{\Phi,S_{\Phi}^{-1}\Phi,\Lambda}(S_{\Phi}^{-1}\phi_i+S_{\Phi}^{-1}\phi_j)&=\langle S_{\Phi}^{-1}\phi_i+S_{\Phi}^{-1}\phi_j,\phi_i\rangle S_{\Phi}^{-1}\phi_i+\langle S_{\Phi}^{-1}\phi_i+S_{\Phi}^{-1}\phi_j,\phi_j\rangle S_{\Phi}^{-1}\phi_j\\[1mm]
    &=\langle S_{\Phi}^{-1}\phi_i,\phi_i\rangle S_{\Phi}^{-1}\phi_i+\langle S_{\Phi}^{-1}\phi_j,\phi_i\rangle S_{\Phi}^{-1}\phi_i+\langle S_{\Phi}^{-1}\phi_i,\phi_j\rangle S_{\Phi}^{-1}\phi_j\\
    &\qquad\qquad\qquad\qquad+\langle S_{\Phi}^{-1}\phi_j,\phi_j\rangle S_{\Phi}^{-1}\phi_j\\[1mm]
     &= \|S_{\Phi}^{-1/2}\phi_i\|^2 S_{\Phi}^{-1}\phi_i+\langle S_{\Phi}^{-1/2}\phi_j,S_{\Phi}^{-1/2}\phi_i\rangle S_{\Phi}^{-1}\phi_i\\[1mm]
&\qquad\qquad+\langle S_{\Phi}^{-1/2}\phi_i,S_{\Phi}^{-1/2}\phi_j\rangle S_{\Phi}^{-1}\phi_j+\|S_{\Phi}^{-1/2}\phi_j\|^2S_{\Phi}^{-1}\phi_j\\[1mm]
      &=\|D_1^{-1/2}D_1M_1^*(e_i)\|^2S_{\Phi}^{-1}\phi_i+\langle D_1^{-1/2}D_1M_1^*(e_j),D_1^{-1/2}D_1M_1^*(e_i)\rangle S_{\Phi}^{-1}\phi_i\\[1mm]
      &\qquad\qquad+\langle D_1^{-1/2}D_1M_1^*(e_i),D_1^{-1/2}D_1M_1^*(e_j)\rangle S_{\Phi}^{-1}\phi_j\\[1mm]
     &\qquad\qquad\qquad  +\|D_1^{-1/2}D_1M_1^*(e_j)\|^2S_{\Phi}^{-1}\phi_j\\[1mm]
      &=\|M_1^*(e_i)\|^2S_{\Phi}^{-1}\phi_i+\langle M_1^*(e_j),M_1^*(e_i)\rangle S_{\Phi}^{-1}\phi_i\\[1mm]
      &\qquad\qquad+\langle M_1^*(e_i),M_1^*(e_j)\rangle S_{\Phi}^{-1}\phi_j+\|M_1^*(e_j)\|^2S_{\Phi}^{-1}\phi_j\\[1mm]
       &=\left(1-\frac{1}{n}\right)S_{\Phi}^{-1}\phi_i-\frac{1}{n}S_{\Phi}^{-1}\phi_i-\frac{1}{n} S_{\Phi}^{-1}\phi_j+\left(1-\frac{1}{n}\right)S_{\Phi}^{-1}\phi_j\\[1mm]
      &=\frac{n-2}{n}\left(S_{\Phi}^{-1}\phi_i+S_{\Phi}^{-1}\phi_j\right).
\end{align*}
Thus, $\frac{n-2}{n}$ is an eigenvalue of $E_{\Phi,S_{\Phi}^{-1}\Phi,\Lambda}$. Similarly, $E_{\Phi,S_{\Phi}^{-1}\Phi,\Lambda}(S_{\Phi}^{-1}\phi_i-S_{\Phi}^{-1}\phi_j)=\left(1-\frac{1}{n}+\frac{1}{n}\right)\break\left(S_{\Phi}^{-1}\phi_i-S_{\Phi}^{-1}\phi_j\right)=S_{\Phi}^{-1}\phi_i-S_{\Phi}^{-1}\phi_j$. Then, the spectrum of $E_{\Phi,S_{\Phi}^{-1}\Phi,\Lambda}$ is $\left\{1,\frac{n-2}{n}\right\}$. Hence, $\rho(E_{\Phi,S_{\Phi}^{-1}\Phi,\Lambda})=1$ for any $\Lambda\subset[n]$ such that $|\Lambda|=2$. Thus, $\rho_{\Phi,S_{\Phi}^{-1}\Phi}^{(2)}=1$. 
\endproof

In the following theorem, we show that the canonical dual frame of a frame generated by a connected graph is the unique spectrally optimal dual frame for $2$-erasures. 
\begin{thm}\label{thm5.2}
Let $G$ be a connected graph on $n$ vertices and $\Phi=\{\phi_i\}_{i\in[n]}$ be an $L_G(n,n-1)$-frame for $\mathbb{C}^{n-1}$ with frame operator $S_{\Phi}$. If $\Psi=\{\psi_i\}_{i\in[n]}$ is a dual frame of $\Phi$ then $\rho_{\Phi,\Psi}^{(2)}\geq 1$. Further, $\rho_{\Phi}^{(2)} = 1$ and $S_{\Phi}^{-1}\Phi$ is the unique $SOD$-frame for $2$-erasures.
\end{thm}
\proof
Let $L$ be the Laplacian matrix of $G$. Since $\Phi$ is an $L_G(n,n-1)$-frame, there exist an orthogonal matrix $M$ and diagonal matrix $D=diag(\lambda_1,\lambda_2,\ldots,\lambda_{n-1},0)$ such that $L=MDM^*$ and $\Phi=D_1M_1^*$ where $M_1$ is formed by considering the first $n-1$ columns of $M$ and $D_1=diag(\sqrt{\lambda_1},\sqrt{\lambda_2},\ldots,\sqrt{\lambda_{n-1}})$. Let $S_{\Phi}$ be the frame operator of $\Phi$ and $\Psi=\{\psi_i\}_{i\in[n]}$ be any dual frame of $\Phi$. Then by Theorem \ref{thm2.5}, there exist $\mu\in\mathbb{C}^{n-1}$ such that $\psi_i=S_{\Phi}^{-1}\phi_i+\mu$ for $i\in[n]$. Let $\Lambda=\{i,j\}\subset[n]$ be arbitrary. Then
\begin{align*}
    E_{\Phi,\Psi,\Lambda}(\psi_i-\psi_j)&=\langle S_{\Phi}^{-1}\phi_i+\mu-S_{\Phi}^{-1}\phi_j-\mu,\phi_i\rangle (S_{\Phi}^{-1}\phi_i+\mu)\\[1mm]
    &\qquad\qquad\qquad+\langle S_{\Phi}^{-1}\phi_i+\mu-S_{\Phi}^{-1}\phi_j-\mu,\phi_j\rangle (S_{\Phi}^{-1}\phi_j+\mu)\\[1mm]
     &= (\|S_{\Phi}^{-1/2}\phi_i\|^2 -\langle S_{\Phi}^{-1/2}\phi_j,S_{\Phi}^{-1/2}\phi_i\rangle) (S_{\Phi}^{-1}\phi_i+\mu)\\[1mm]
&\qquad\qquad+(\langle S_{\Phi}^{-1/2}\phi_i,S_{\Phi}^{-1/2}\phi_j\rangle -\|S_{\Phi}^{-1/2}\phi_j\|^2)(S_{\Phi}^{-1}\phi_j+\mu)\\[1mm]
      &=(\|M_1^*(e_i)\|^2-\langle M_1^*(e_j),M_1^*(e_i)\rangle)(S_{\Phi}^{-1}\phi_i+\mu)\\[1mm]
      &\qquad\qquad+(\langle M_1^*(e_i),M_1^*(e_j)\rangle -\|M_1^*(e_j)\|^2)(S_{\Phi}^{-1}\phi_j+\mu)\\[1mm]
       &=\left(1-\frac{1}{n}+\frac{1}{n}\right)(S_{\Phi}^{-1}\phi_i+\mu)+\left(-\frac{1}{n}-1+\frac{1}{n}\right)(S_{\Phi}^{-1}\phi_j+\mu)\\[1mm]
      &=\psi_i-\psi_j.
\end{align*}
Thus, $1$ is an eigenvalue of $E_{\Phi,\Psi,\Lambda}$. Then $\rho(E_{\Phi,\Psi,\Lambda})\geq 1$. Hence, $\rho^{(2)}_{\Phi,\Psi}=\max\{\rho(E_{\Phi,\Psi,\Lambda}):\Lambda\subset[n] \text{ and } |\Lambda|=2 \}\geq 1$. Thus, by Theorem \ref{thm4.1} and Theorem \ref{thm5.1}, $S_{\Phi}^{-1}\Phi$ is the unique $SOD$-frame for $2$-erasures.
\endproof

\begin{corollary}
    Let $G$ be a connected graph on $n$ vertices and $\Phi=\{\phi_i\}_{i\in[n]}$ be an $G(n,n-1)$-frame for $\mathbb{C}^{n-1}$ with frame operator $S_{\Phi}$. Then canonical dual of $\Phi$ is the unique $SOD$-frame for $2$-erasures and $\rho^{(2)}_{\Phi}=2$.
\end{corollary}

In the next theorem, we show that if $\Phi$ is an $L_G(n,k)$-frame generated by any graph $G$ (need not be connected) with frame operator $S_{\Phi}$ then $\rho_{\Phi,S_{\Phi}^{-1}\Phi}^{(2)}=1$. 
\begin{thm}\label{thm5.3}
Let $G$ be a simple graph on $n$ vertices with $k$ components and $\Phi=\{\phi_i\}_{i\in[n]}$ is an $L_G(n,n-k)$-frame for $\mathbb{C}^k$. If $S_{\Phi}$ is the frame operator of $\Phi$ then $\rho_{\Phi,S_{\Phi}^{-1}\Phi}^{(2)}=1$. 
\end{thm}
\proof

	Let $G_1,G_2,\ldots,G_k$ be the components of $G$ and $ |V(G_i)|=n_i $, for $ i\in [k] $. Suppose $V(G_i)=\{l_{i-1}+1,l_{i-1}+2,\ldots,l_i\}$ where $ l_j=n_1+\cdots+n_{j} $ for $ j\in [k]$ and $l_0=0$. For $i\in [k]$, let $\lambda_{i,1}\geq\lambda_{i,2}\geq\cdots\geq\lambda_{i,n_i-1}>\lambda_{i,n_i}=0$  be the eigenvalues of $L(G_i)$. Then there exist an orthogonal matricx $M_i$ such that $L(G_i)=M_iD_iM_i^*$ where $D_i=diag(\lambda_{i,1},\lambda_{i,2},\ldots,\lambda_{i, n_i-1},0)$, for $i\in[k]$, and $L(G)=MDM^*$ where $M=M_1\oplus\cdots\oplus M_k$ and $D=D_1\oplus\cdots\oplus D_k$. Thus, the synthesis operator of an $L_G(n,n-k)$-frame, say $\widetilde{\Phi}=\{ \widetilde{\phi}_i\}_{i\in[n]}$ is $B=B_1\oplus B_2\oplus\cdots\oplus B_k$ such that  $B_j=diag(\sqrt{\lambda_{j,1}},\sqrt{\lambda_{j,2}},\ldots,\sqrt{\lambda_{j,n_j-1}})\widetilde{M_j}^*$ (where $\widetilde{M_j}^*$ is formed by taking the first $n_j-1$ columns of $M_j$) is the synthesis operator of the $L_{G_i}(n_i,n_i-1)$-frame, say $\Phi_j=\{\phi_{j,i}\}_{i\in[n_j]}$ for $\mathbb{C}^{n_j-1}$, for $j\in[k]$. By Theorem \ref{thm2.3}, there is a unitary operator $ V $ such that $ \phi_i=V\widetilde{\phi}_i $ for all $i\in[n]$.  Let $S_{\Phi}$ and $S_{\widetilde{\Phi}}$  be the frame operators of $\Phi$  and $\widetilde{\Phi}$, respectively. For $i\in [k]$, let $S_{\Phi_i}$ be the frame operator of $\Phi_i$. Then
	\begin{align*}
		S_{\widetilde{\Phi}}&=BB^*\\
		&=B_1B_1^*\oplus\cdots\oplus B_kB_k^*\\
		&=S_{\Phi_1}\oplus\cdots\oplus S_{\Phi_k}.	
	\end{align*}
	Thus, $S_{\widetilde{\Phi}}^{-1}=S_{\Phi_1}^{-1}\oplus\cdots\oplus S_{\Phi_k}^{-1}$, and $ S_{\Phi}^{-1}=VS_{\widetilde{\Phi}}^{-1}V^* $. Let $\Lambda=\{p,q\}\subset [n]$ be arbitrary. Now we proceed in the following cases.

\noindent\textbf{Case I: } Let $\Lambda\subset V(G_j)$ for some $j\in[k]$. Suppose $f=\left[\begin{array}{c}
	f_1 \\
f_2\\
\vdots\\
f_k
	\end{array}\right]\in\mathbb{C}^{n-k}$ where $f_i\in\mathbb{C}^{n_i-1}$. We have
 \begin{align*}
     E_{\widetilde{\Phi},S_{\widetilde{\Phi}}^{-1}\widetilde{\Phi},\Lambda}(f)&=\langle f,\widetilde{\phi}_p\rangle S_{\widetilde{\Phi}}^{-1}\widetilde{\phi}_p+\langle f,\widetilde{\phi}_q\rangle S_{\widetilde{\Phi}}^{-1}\widetilde{\phi}_q\\
&=\langle f_j,\widetilde{\phi}_{j,p-l_{j-1}}\rangle S_{\widetilde{\Phi}}^{-1}\widetilde{\phi}_p+\langle f_j,\widetilde{\phi}_{j,q-l_{j-1}}\rangle S_{\widetilde{\Phi}}^{-1}\widetilde{\phi}_q\\
&=\left[\begin{array}{c}
   0 \\
   \vdots\\
   0\\
   E_{\Phi_j,S_{\Phi_j}^{-1}\Phi_j,\Lambda'}(f_j)\\
   0\\
   \vdots\\
   0
\end{array}
\right], \text{ where } \Lambda'=\{p-l_{j-1},q-l_{j-1}\}.
 \end{align*}
Thus, by Theorem \ref{thm5.1}, $\rho(E_{\widetilde{\Phi},S_{\widetilde{\Phi}}^{-1}\widetilde{\Phi},\Lambda})=\rho(E_{\Phi_j,S_{\Phi_j}^{-1}\Phi_j,\Lambda'})=1$.
 
\noindent\textbf{Case II: } Let $p\in V(G_s)$ and $q\in V(G_r)$ where $s\neq r$. Then, we have
\begin{align*}
     E_{\widetilde{\Phi},S_{\widetilde{\Phi}}^{-1}\widetilde{\Phi},\Lambda}(S_{\widetilde{\Phi}}^{-1}\widetilde{\phi}_p)&=\langle S_{\widetilde{\Phi}}^{-1}\widetilde{\phi}_p,\widetilde{\phi}_p\rangle S_{\widetilde{\Phi}}^{-1}\widetilde{\phi}_p+\langle S_{\widetilde{\Phi}}^{-1}\widetilde{\phi}_p,\widetilde{\phi}_q\rangle S_{\widetilde{\Phi}}^{-1}\widetilde{\phi}_q\\[1mm]
&=\|S_{\widetilde{\Phi}}^{-1/2}\widetilde{\phi}_p\|^2S_{\widetilde{\Phi}}^{-1}\widetilde{\phi}_p+0(S_{\widetilde{\Phi}}^{-1}\widetilde{\phi}_q)\\[1mm]
&=\frac{n_s-1}{n_s}S_{\widetilde{\Phi}}^{-1}\widetilde{\phi}_p.
 \end{align*}
 Similarly, $E_{\widetilde{\Phi},S_{\widetilde{\Phi}}^{-1}\widetilde{\Phi},\Lambda}(S_{\widetilde{\Phi}}^{-1}\widetilde{\phi}_q)=\frac{n_r-1}{n_r}S_{\widetilde{\Phi}}^{-1}\widetilde{\phi}_q$. Thus, $\rho(E_{\widetilde{\Phi},S_{\widetilde{\Phi}}^{-1}\widetilde{\Phi},\Lambda})=\max\{\frac{n_s-1}{n_s},\frac{n_r-1}{n_r}\}<1$. 
 
 Hence, by Case I and Case II, $\max\{\rho(E_{\widetilde{\Phi},S_{\widetilde{\Phi}}^{-1}\widetilde{\Phi},\Lambda}):\Lambda\subset[n]\text{ and }|\Lambda|=2\}=1$. Therefore, we have
 \begin{align*}
     \rho^{(2)}_{\Phi,S_{\Phi}^{-1}\Phi}&=\max\{\rho(E_{\Phi,S_{\Phi}^{-1}\Phi,\Lambda}):\Lambda\subset[n]\text{ and }|\Lambda|=2\}\\
     &=\max\{\rho(VE_{\widetilde{\Phi},S_{\widetilde{\Phi}}^{-1}\widetilde{\Phi},\Lambda}V^*):\Lambda\subset[n]\text{ and }|\Lambda|=2\}\\
      &=\max\{\rho(E_{\widetilde{\Phi},S_{\widetilde{\Phi}}^{-1}\widetilde{\Phi},\Lambda}):\Lambda\subset[n]\text{ and }|\Lambda|=2\}=1.
 \end{align*}
\endproof
In Theorem \ref{thm4.3}, we show that the canonical dual frames of a frame generated by disconnected graphs are non-unique $SOD$-frames for $1$-erasure. Hence, it would be interesting to see whether canonical dual frames of frames generated by disconnected graphs are $SOD$-frames for $2$-erasures or not. In the next two theorems, we show that canonical duals of $L_G(n,k)$-frames, where $G$ is a disconnected graph, are non-unique $SOD$-frames for $2$-erasures.
\begin{thm}\label{thm5.4}
	Suppose $G$ is a graph with $n$ vertices and $k>1$ components. If $\Phi=\{\phi_i\}_{i\in[n]}$ is an $L_G(n,n-k)$-frame for $\mathbb{C}^{n-k}$ with frame operator $S_{\Phi}$, then the canonical dual frame $S_{\Phi}^{-1}\Phi$ of $\Phi$ is a spectrally optimal dual frame of $\Phi$ for $2$-erasures. Further, $\rho^{(2)}_{\Phi}=1$.
\end{thm}	
\proof
Let $G_1,G_2,\ldots,G_k$ be the connected components of $G$ such that $|V(G_i)|=n_i$ for $i\in[k]$. Let the vertex set of $G_i$ be $\{l_{i-1}+1,l_{i-1}+2,\ldots,l_i\}$ for $i\in[k]$ where $l_0=0$ and for $1\leq i\leq k$, $l_i=n_1+n_2+\cdots+n_i$. Let $\Psi=\{\psi_i\}_{i\in[n]}$ be any dual frame of $\Phi$. Then by Theorem \ref{thm2.5}, there exist $\mu_1,\mu_2,\ldots,\mu_k\in\mathbb{C}^{n-k}$ such that $\Psi=\bigcup\limits_{j\in[k]}\{S_{\Phi}^{-1}\phi_{l_{j-1}+i}+\mu_j\}_{i\in[n_j]}$. Without loss of generality, assume that $|V(G_1)|\geq 2$. Take $\Lambda=\{1,2\}$. Then, $E_{\Phi,\Psi, \Lambda}(f)=\langle f, \phi_1\rangle 
 (S_{\Phi}^{-1}\phi_1+\mu_1)+ \langle f, \phi_2\rangle (S_{\Phi}^{-1}\phi_2+\mu_1)$. Thus, we have
\begin{align*}
E_{\Phi,\Psi, \Lambda}(S_{\Phi}^{-1}\phi_1+\mu_1-(S_{\Phi}^{-1}\phi_2+\mu_1))&=\langle S_{\Phi}^{-1}\phi_1+\mu_1-S_{\Phi}^{-1}\phi_2-\mu_1, \phi_1\rangle (S_{\Phi}^{-1}\phi_1+\mu_1)\\[1mm]
&\qquad\qquad+ \langle S_{\Phi}^{-1}\phi_1+\mu_1-S_{\Phi}^{-1}\phi_2-\mu_1, \phi_2\rangle (S_{\Phi}^{-1}\phi_2+\mu_1)\\[1mm]
&=(\langle S_{\Phi}^{-1}\phi_1, \phi_1\rangle-\langle S_{\Phi}^{-1}\phi_2, \phi_1\rangle) (S_{\Phi}^{-1}\phi_1+\mu_1)\\[1mm]
&\qquad\qquad+ (\langle S_{\Phi}^{-1}\phi_1, \phi_2\rangle-\langle S_{\Phi}^{-1}\phi_2, \phi_2\rangle) (S_{\Phi}^{-1}\phi_2+\mu_1)\\[1mm]
&=(\|S_{\Phi}^{-1/2}\phi_1\|^2-\langle S_{\Phi}^{-1/2}\phi_2, S_{\Phi}^{-1/2}\phi_1\rangle) (S_{\Phi}^{-1}\phi_1+\mu_1)\\[1mm]
&\qquad\qquad+ (\langle S_{\Phi}^{-1/2}\phi_1, S_{\Phi}^{-1/2}\phi_2\rangle-\|S_{\Phi}^{-1/2}\phi_2\|^2) (S_{\Phi}^{-1}\phi_2+\mu_1)\\[1mm]
&=\left(1-\frac{1}{n_1}+\frac{1}{n_1}\right)(S_{\Phi}^{-1}\phi_1+\mu_1)-\left(1-\frac{1}{n_1}+\frac{1}{n_1}\right)(S_{\Phi}^{-1}\phi_2+\mu_1)\\[1mm]
&=(S_{\Phi}^{-1}\phi_1+\mu_1)-(S_{\Phi}^{-1}\phi_2+\mu_1).
\end{align*}
Thus, $1$ is an eigenvalue of $E_{\Phi,\Psi, \Lambda}$. Hence, $\rho(E_{\Phi,\Psi, \Lambda})\geq 1$. This gives that $\rho^{(2)}_{\Phi,\Psi}\geq\rho(E_{\Phi,\Psi, \Lambda})\geq 1$. By Theorem \ref{thm4.3}, $S_{\Phi}^{-1}\Phi$ is spectrally optimal dual frame for $1$-erasures and by Theorem \ref{thm5.3}, $\rho^{(2)}_{\Phi,S_{\Phi}^{-1}\Phi}=1\geq \rho^{(2)}_{\Phi,\Psi}$ for any dual frame $\Psi$ of $\Phi$. Thus, the canonical dual frame of $\Phi$ is a spectrally optimal dual frame for $2$-erasures.
\endproof

In the next theorem, we show the existence of an alternate dual frame of a frame generated by a disconnected graph such that the alternate dual is a spectrally optimal dual frame for $2$-erasures.
\begin{thm}\label{thm5.5}
	Suppose $G$ is a graph with $n$ vertices and $k>1$ connected components. If $\Phi=\{\phi_i\}_{i\in[n]}$ is an $L_G(n,n-k)$-frame for $\mathbb{C}^{n-k}$ with frame operator $S_{\Phi}$, then there exist an alternate dual frame of $\Phi$ such that it is a spectrally optimal dual frame of $\Phi$ for $2$-erasures.
\end{thm}	
\proof
Let $G_1,G_2,\ldots,G_k$ be the connected components of $G$ such that $|V(G_i)|=n_i$ for $i\in[k]$. Let the vertex set of $G_i$ be $\{v_{l_{i-1}+1},v_{l_{i-1}+2},\ldots,v_{l_i}\}$ for $i\in[k]$ where $l_0=0$ and for $i\geq 1$, $l_i=n_1+n_2+\cdots+n_i$. Let $\widetilde{\Phi}$ be the $L_G(n,k)$-frame as obtained in Theorem \ref{thm5.3}. By Theorem \ref{thm5.4}, $S_{\widetilde{\Phi}}^{-1}\widetilde{\Phi}$ is a spectrally optimal dual frame of $\widetilde{\Phi}$ for $2$-erasures and $\rho^{(2)}_{\widetilde{\Phi},S_{\widetilde{\Phi}}^{-1}\widetilde{\Phi}}=1$. Let $\mu=[\mu_1\, , \mu_2\, ,\ldots,\, \mu_{n-k}]^t\in\mathbb{C}^{n-k}$ such that $\mu_i=\begin{cases}
    1,\, i=l_1+1\\
    0, \, \text{otherwise}
\end{cases}$. Then $\Psi=\{\psi_i\}_{i\in[n]}=\{S_{\widetilde{\Phi}}^{-1}\widetilde{\phi}_i+\mu\}_{i\in[l_1]}\bigcup\{S_{\widetilde{\Phi}}^{-1}\widetilde{\phi}_i\}_{i\in[n]\setminus[l_1]}$ is an alternate dual frame of $\widetilde{\Phi}$. Then 
\begin{align*}
\rho^{(1)}_{\widetilde{\Phi},\Psi}&=\max\{|\langle \psi_i,\widetilde{\phi}_i\rangle|:i\in[n]\}\\
&=\max\left(\{|\langle S_{\widetilde{\Phi}}^{-1}\widetilde{\phi}_i+\mu,\widetilde{\phi}_i\rangle|:i\in[l_1]\}\cup\{|\langle S_{\Phi}^{-1}\phi_i,\widetilde{\phi}_i\rangle|:i\in[n]\setminus[l_1]\}\right)\\
&=\max\left(\{|\langle S_{\Phi}^{-1}\phi_i,\widetilde{\phi}_i\rangle+\langle \mu,\widetilde{\phi}_i\rangle|:i\in[l_1]\}\cup\{|\langle S_{\Phi}^{-1}\phi_i,\widetilde{\phi}_i\rangle|:i\in[n]\setminus[l_1]\}\right) \\
&=\max\left(\{|\langle S_{\Phi}^{-1}\phi_i,\widetilde{\phi}_i\rangle|:i\in[l_1]\}\cup\{|\langle S_{\Phi}^{-1}\phi_i,\widetilde{\phi}_i\rangle|:i\in[n]\setminus[l_1]\}\right)=\rho^{(1)}_{\widetilde{\Phi},S^{-1}_{\widetilde{\Phi}}\widetilde{\Phi}}.
\end{align*}
Thus, by Theorem \ref{thm4.3}, $\Psi$ is a spectrally optimal dual frame of $\widetilde{\Phi}$ for $1$-erasure. Let $\Lambda=\{p,q\}\subset[n]$ be arbitrary. 
Then $E_{\widetilde{\Phi},\Psi, \Lambda}(f)=\langle f,\widetilde{\phi}_p\rangle \psi_p+\langle f,\widetilde{\phi}_q\rangle \psi_q$. Now we will proceed in cases. 

\noindent\textbf{Case I: } Let $\Lambda\subseteq V(G_j)$ for some $j\in[k]$. Then, we have
\begin{align*}
    E_{\widetilde{\Phi},\Psi, \Lambda}(\psi_p+\psi_q)&=\langle \psi_p+\psi_q,\widetilde{\phi}_p\rangle \psi_p+\langle \psi_p+\psi_q,\widetilde{\phi}_q\rangle \psi_q\\[2mm]
   &=\begin{cases}
       \langle S_{\widetilde{\Phi}}^{-1}\widetilde{\phi}_p+ S_{\widetilde{\Phi}}^{-1}\widetilde{\phi}_q+2\mu,\widetilde{\phi}_p\rangle \psi_p+\langle S_{\widetilde{\Phi}}^{-1}\widetilde{\phi}_p+ S_{\widetilde{\Phi}}^{-1}\widetilde{\phi}_p+2\mu,\widetilde{\phi}_q\rangle \psi_q \text{ if } i=1\\
        \langle S_{\widetilde{\Phi}}^{-1}\widetilde{\phi}_p+ S_{\widetilde{\Phi}}^{-1}\widetilde{\phi}_q,\widetilde{\phi}_p\rangle \psi_p+\langle S_{\widetilde{\Phi}}^{-1}\widetilde{\phi}_p+ S_{\widetilde{\Phi}}^{-1}\widetilde{\phi}_q,\widetilde{\phi}_q\rangle \psi_q \text{ otherwise } 
   \end{cases}\\[2mm]
   &=\begin{cases}
       (\|S_{\widetilde{\Phi}}^{-1/2}\widetilde{\phi}_p\|^2+\langle S_{\widetilde{\Phi}}^{-1/2}\widetilde{\phi}_q,S_{\widetilde{\Phi}}^{-1/2}\widetilde{\phi}_p\rangle+2\langle \mu,\widetilde{\phi}_p\rangle) \psi_p
       \\\qquad\qquad+(\|S_{\widetilde{\Phi}}^{-1/2}\widetilde{\phi}_q\|^2+\langle S_{\widetilde{\Phi}}^{-1/2}\widetilde{\phi}_p,S_{\widetilde{\Phi}}^{-1/2}\widetilde{\phi}_q\rangle+2\langle \mu,\widetilde{\phi}_q\rangle) \psi_q \text{ if } i=1\\
        (\|S_{\widetilde{\Phi}}^{-1/2}\widetilde{\phi}_p\|^2+\langle S_{\widetilde{\Phi}}^{-1/2}\widetilde{\phi}_q,S_{\widetilde{\Phi}}^{-1/2}\widetilde{\phi}_p\rangle) \psi_p
       \\\qquad\qquad+(\|S_{\widetilde{\Phi}}^{-1/2}\widetilde{\phi}_q\|^2+\langle S_{\widetilde{\Phi}}^{-1/2}\widetilde{\phi}_p,S_{\widetilde{\Phi}}^{-1/2}\widetilde{\phi}_q\rangle) \psi_q  \text{ otherwise } 
   \end{cases}\\[2mm]
   &=\begin{cases}
       \left(\left(1-\frac{1}{n_1}\right)-\frac{1}{n_1}+0\right) \psi_p+\left(\left(1-\frac{1}{n_1}\right)-\frac{1}{n_1}+0\right) \psi_q \text{ if } j=1\\
       \left(\left(1-\frac{1}{n_j}\right)-\frac{1}{n_j}\right) \psi_p+\left(\left(1-\frac{1}{n_j}\right)-\frac{1}{n_j}\right) \psi_q  \text{ otherwise } 
   \end{cases}\\[2mm]
   &= \left(1-\frac{2}{n_j}\right)(\psi_p+\psi_q).
\end{align*}
Similarly, $E_{\widetilde{\Phi},\Psi, \Lambda}(\psi_p-\psi_q)=\psi_p-\psi_q$. Thus, $\rho(E_{\widetilde{\Phi},\Psi, \Lambda})=1$.

\noindent\textbf{Case II: } Let $p\in V(G_r)$ and $q\in V(G_s)$ where $r,s\in[k]\setminus\{1\}$ and $r\neq s$. Then, we have
\begin{align*}
    E_{\widetilde{\Phi},\Psi, \Lambda}(\psi_p)&=\langle \psi_p,\widetilde{\phi}_p\rangle \psi_p+\langle \psi_p,\widetilde{\phi}_q\rangle \psi_q\\[1mm]
    &=\langle S_{\widetilde{\Phi}}^{-1}\widetilde{\phi}_p,\widetilde{\phi}_p\rangle \psi_p+\langle S_{\widetilde{\Phi}}^{-1}\widetilde{\phi}_p,\widetilde{\phi}_q\rangle \psi_q\\[1mm]
&=\|S_{\widetilde{\Phi}}^{-1/2}\widetilde{\phi}_p\|^2\psi_p+\langle S_{\widetilde{\Phi}}^{-1/2}\widetilde{\phi}_p,S_{\widetilde{\Phi}}^{-1/2}\widetilde{\phi}_q\rangle \psi_q\\[1mm]
&=\left(1-\frac{1}{n_r}\right)\psi_p+0\\[1mm]
&=\left(1-\frac{1}{n_r}\right)\psi_p.
\end{align*}
Similarly, $E_{\widetilde{\Phi},\Psi, \Lambda}(\psi_q)=\left(1-\frac{1}{n_s}\right)\psi_q$. Hence, $\rho(E_{\widetilde{\Phi},\Psi, \Lambda})=\max\left\{1-\frac{1}{n_r},1-\frac{1}{n_s}\right\}<1$.

\noindent\textbf{Case III: } Let $p\in V(G_1)$ and $q\in V(G_r)$ where $r\neq 1$. Then $E_{\widetilde{\Phi},\Psi, \Lambda}(f)=\langle f, \widetilde{\phi}_p\rangle \psi_p+\langle f, \widetilde{\phi}_q\rangle \psi_q=\langle f, \widetilde{\phi}_p\rangle (S_{\widetilde{\Phi}}^{-1}\widetilde{\phi}_p+\mu)+\langle f, \widetilde{\phi}_q\rangle S_{\widetilde{\Phi}}^{-1}\widetilde{\phi}_q$. The matrix representation of $E_{\widetilde{\Phi},S_{\widetilde{\Phi}}^{-1}\widetilde{\Phi}, \Lambda}$ with respect to the standard canonical orthonormal basis is of the form, $[E_{\widetilde{\Phi},S_{\widetilde{\Phi}}^{-1}\widetilde{\Phi}, \Lambda}]=\left[\begin{array}{cc}
    E_1 &0 \\
0 & E_2
\end{array}\right]$ for some matrices $E_1$ of order $n_1-1$ and $E_2$ of order $(n-k)-(n_1-1)$. Then, the matrix representation of $E_{\widetilde{\Phi},\Psi, \Lambda}$ with respect to the standard canonical orthonormal basis is $[E_{\widetilde{\Phi},\Psi, \Lambda}]=\left[\begin{array}{cc}
    E_1 &0 \\
E_3 & E_2
\end{array}\right]$ where $E_3=\left[\begin{array}{cccc}
   \alpha_1  &  \alpha_2 & \cdots & \alpha_{l_1-1}\\
  0  &  0 & \cdots & 0\\
   \vdots  &  \vdots & \ddots & \vdots\\
    0  &  0 & \cdots & 0
\end{array}\right]$ and $\alpha_i$ is the $i^{th}$-component of $\widetilde{\phi}_p$. Thus, $\rho(E_{\widetilde{\Phi},\Psi, \Lambda})=\rho([E_{\widetilde{\Phi},S_{\widetilde{\Phi}}^{-1}\widetilde{\Phi}, \Lambda}])=1$. Then, by Theorem \ref{thm5.4}, $\rho^{(2)}_{\widetilde{\Phi},\Psi}=1=\rho^{(2)}_{\widetilde{\Phi}}$. Hence, $\Psi$ is an alternate dual frame of $\widetilde{\Phi}$ such that $\Psi$ is a spectrally optimal dual frame of $\widetilde{\Phi}$ for $2$-erasures. Therefore, by Theorem \ref{thm4.2}, there exists an alternate dual frame of $\Phi$, which is an $SOD$-frame of $\Phi$ for $2$-erasures.
\endproof

\begin{corollary}
    Suppose $G$ is a graph with $n$ vertices and $k>1$ connected components. If $\Phi=\{\phi_i\}_{i\in[n]}$ is an $G(n,n-k)$-frame for $\mathbb{C}^{n-k}$ with frame operator $S_{\Phi}$, then there exist an alternate dual frame of $\Phi$ such that it is an $SOD$-frame of $\Phi$ for $2$-erasures.
\end{corollary}
	In the following example, we show explicitly the existence of an alternate dual frame that is an $SOD$-frame for $2$-erasures of a frame generated by a disconnected graph.
\begin{example}\label{Exa2}{\em 

    Consider the frame $Phi=\{\phi_i\}_{i\in[5]}$ generated by graph $G$ given in Example \ref{Exa1}. The canonical dual frame $S_{\Phi}^{-1}\Phi$ is an $SOD$-frame of $\Phi$ for $1$-erasure. Now, we will show that $S_{\Phi}^{-1}\Phi$ is a non-unique $SOD$-frame of $\Phi$ for $2$-erasures. First, we give the matrix representation of error operators. For $\Lambda=\{i,j\}\subset[5]$, $E_{\Phi,S_{\Phi}^{-1}\Phi,\Lambda}(f)=\langle f,\phi_i\rangle S_{\Phi}^{-1}\phi_i+\langle f,\phi_j\rangle S_{\Phi}^{-1}\phi_j$. We have
     $$\left[E_{\Phi,S_{\Phi}^{-1}\Phi,\{1,2\}}\right]=\left[\begin{array}{ccc}
    \frac{1}{2} & \frac{-1}{2\sqrt{3}} & 0\\[1mm]
    \frac{-1}{2\sqrt{3}} & \frac{5}{6} &  0\\[1mm]
    0 & 0 & 0
    \end{array}\right], \left[E_{\Phi,S_{\Phi}^{-1}\Phi,\{1,3\}}\right]=\left[\begin{array}{ccc}
    1 & 0 & 0\\[1mm]
    0 & \frac{1}{3} &  0\\[1mm]
    0 & 0 & 0
    \end{array}\right],
    \left[E_{\Phi,S_{\Phi}^{-1}\Phi,\{2,3\}}\right]=\left[\begin{array}{ccc}
    \frac{1}{2} & \frac{1}{2\sqrt{3}} & 0\\[1mm]
    \frac{1}{2\sqrt{3}} & \frac{5}{6} &  0\\[1mm]
    0 & 0 & 0
    \end{array}\right],$$ $$\left[E_{\Phi,S_{\Phi}^{-1}\Phi,\{4,5\}}\right]=\left[\begin{array}{ccc}
    0 & 0 & 0\\[1mm]
    0 & 0 &  0\\[1mm]
    0 & 0 & 1
    \end{array}\right], \left[E_{\Phi,S_{\Phi}^{-1}\Phi,\{1,4\}}\right]=\left[E_{\Phi,S_{\Phi}^{-1}\Phi,\{1,5\}}\right]=\left[\begin{array}{ccc}
    \frac{1}{2} & \frac{-1}{2\sqrt{3}} & 0\\[1mm]
    \frac{-1}{2\sqrt{3}} & \frac{1}{6} &  0\\[1mm]
    0 & 0 & \frac{1}{2}
    \end{array}\right],$$ $$\left[E_{\Phi,S_{\Phi}^{-1}\Phi,\{2,4\}}\right]=\left[E_{\Phi,S_{\Phi}^{-1}\Phi,\{2,5\}}\right]=\left[\begin{array}{ccc}
    0 & 0 & 0\\[1mm]
    0 & \frac{2}{3} &  0\\[1mm]
    0 & 0 & \frac{1}{2}
    \end{array}\right], \left[E_{\Phi,S_{\Phi}^{-1}\Phi,\{3,4\}}\right]=\left[E_{\Phi,S_{\Phi}^{-1}\Phi,\{3,5\}}\right]=\left[\begin{array}{ccc}
    \frac{1}{2} & \frac{1}{2\sqrt{3}} & 0\\[1mm]
    \frac{1}{2\sqrt{3}} & \frac{1}{6} &  0\\[1mm]
    0 & 0 & \frac{1}{2}
    \end{array}\right].$$ Now $\rho(E_{\Phi,S_{\Phi}^{-1}\Phi,\{1,2\}})=\rho(E_{\Phi,S_{\Phi}^{-1}\Phi,\{1,3\}})=\rho(E_{\Phi,S_{\Phi}^{-1}\Phi,\{2,3\}})=\max\left\{1,\frac{1}{3},0\right\}=1$, $\rho(E_{\Phi,S_{\Phi}^{-1}\Phi,\{4,5\}})=\max\left\{1,0\right\}=1$ and $\rho(E_{\Phi,S_{\Phi}^{-1}\Phi,\{1,4\}})=\rho(E_{\Phi,S_{\Phi}^{-1}\Phi,\{1,5\}})=\rho(E_{\Phi,S_{\Phi}^{-1}\Phi,\{2,4\}})=\rho(E_{\Phi,S_{\Phi}^{-1}\Phi,\{2,5\}})=\rho(E_{\Phi,S_{\Phi}^{-1}\Phi,\{3,4\}})=\rho(E_{\Phi,S_{\Phi}^{-1}\Phi,\{3,5\}})=\max\left\{\frac{2}{3},\frac{1}{3},0\right\}=\frac{2}{3}$. Thus, $\rho_{\Phi,S_{\Phi}^{-1}\Phi}^{(2)}=\max\left\{1,\frac{2}{3}\right\}=1$. Let $\Psi=\{\psi_i\}_{i\in[5]}$ be any dual frame of $\Phi$. Then there exist scalars $a_o,b_o,c_o,e_o,f_o,g_o$ such that $\Psi=\{\psi_i\}_{i\in[5]}=\left\{\left[\begin{array}{c}
    \frac{1}{\sqrt{6}} + a_o\\[1mm]
    \frac{-1}{3\sqrt{2}} + b_o \\[1mm]
    c_o
    \end{array}\right], \left[\begin{array}{c}
    a_o \\[1mm]
    \frac{2}{3\sqrt{2}} + b_o\\[1mm]
    c_o 
    \end{array}\right], \left[\begin{array}{c}
    \frac{-1}{\sqrt{6}} + a_o \\[1mm]
    \frac{-1}{3\sqrt{2}} + b_o \\[1mm]
    c_o 
    \end{array}\right], \left[\begin{array}{c}
    e_o\\[1mm]
    f_o\\[1mm]
    \frac{1}{2} +g_o
    \end{array}\right], \left[\begin{array}{c}
    e_o\\[1mm]
    f_o\\[1mm]
    \frac{-1}{2} + g_o
    \end{array}\right]\right\}$. \\[2mm] For $\Lambda=\{1,2\}$, $E_{\Phi,\Psi,\Lambda}(f)=\langle f,\phi_1\rangle \psi_1+\langle f,\phi_2\rangle \psi_2$. For $f=\left[\begin{array}{c}
    \frac{1}{\sqrt{6}}\\[1mm]
    \frac{-1}{\sqrt{2}}\\[1mm]
    0
    \end{array}\right]$, $E_{\Phi,\Psi,\Lambda}(f)=\left[\begin{array}{c}
    \frac{1}{\sqrt{6}}+a_o\\[1mm]
    \frac{-1}{3\sqrt{2}}+b_o\\[1mm]
    c_o
    \end{array}\right]-\left[\begin{array}{c}
    a_o\\[1mm]
    \frac{2}{3\sqrt{2}}+b_o\\[1mm]
    c_o
    \end{array}\right]=\left[\begin{array}{c}
    \frac{1}{\sqrt{6}}\\[1mm]
    \frac{-1}{\sqrt{2}}\\[1mm]
    0
    \end{array}\right]$.
Thus, $1$ is an eigenvalue of $E_{\Phi,\Psi,\Lambda}$ and hence $\rho(E_{\Phi,\Psi,\Lambda})\geq 1$. Therefore, $\rho^{(2)}_{\Phi,\Psi}\geq 1=\rho^{(2)}_{\Phi,S^{-1}_{\Phi}\Phi}$ for any dual frame $\Psi$ of the frame $\Phi$. Thus, $S^{-1}_{\Phi}\Phi$ is an $SOD$-frame of $\Phi$ for $2$-erasures. Consider the dual frame $\Psi_1$ given in Example \ref{Exa1}. It is shown in Example \ref{Exa1} that $\Psi_1$ is an $SOD$-frame of $\Phi$ for $1$-erasure. Further, we have 
    $$\left[E_{\Phi,\Psi_1,\{1,2\}}\right]=\left[\begin{array}{ccc}
    \frac{1}{2} & \frac{-1}{2\sqrt{3}} & 0\\
    \frac{-1}{2\sqrt{3}} & \frac{5}{6} &  0\\
    \frac{\sqrt{3}}{\sqrt{2}} & \frac{1}{\sqrt{2}} & 0
    \end{array}\right], \left[E_{\Phi,\Psi_1,\{1,3\}}\right]=\left[\begin{array}{ccc}
    1 & 0 & 0\\
    0 & \frac{1}{3} &  0\\
    0 & -\sqrt{2} & 0
    \end{array}\right],
    \left[E_{\Phi,\Psi_1,\{2,3\}}\right]=\left[\begin{array}{ccc}
    \frac{1}{2} & \frac{1}{2\sqrt{3}} & 0\\
    \frac{1}{2\sqrt{3}} & \frac{5}{6} &  0\\
    \frac{-\sqrt{3}}{\sqrt{2}} & \frac{1}{\sqrt{2}} & 0
    \end{array}\right],$$ $$\left[E_{\Phi,\Psi_1,\{4,5\}}\right]=\left[\begin{array}{ccc}
    0 & 0 & 0\\
    0 & 0 &  0\\
    0 & 0 & 1
    \end{array}\right], \left[E_{\Phi,\Psi_1,\{1,4\}}\right]=\left[E_{\Phi,\Psi_1,\{1,5\}}\right]=\left[\begin{array}{ccc}
    \frac{1}{2} & \frac{-1}{2\sqrt{3}} & 0\\
    \frac{-1}{2\sqrt{3}} & \frac{1}{6} &  0\\
    \frac{\sqrt{3}}{\sqrt{2}} & \frac{1}{\sqrt{2}}& \frac{1}{2}
    \end{array}\right],$$ $$\left[E_{\Phi,\Psi_1,\{2,4\}}\right]=\left[E_{\Phi,\Psi_1,\{2,5\}}\right]=\left[\begin{array}{ccc}
    0 & 0 & 0\\
    0 & \frac{2}{3} &  0\\
    0 & \sqrt{2} & \frac{1}{2}
    \end{array}\right], \left[E_{\Phi,\Psi_1,\{3,4\}}\right]=\left[E_{\Phi,\Psi_1,\{3,5\}}\right]=\left[\begin{array}{ccc}
    \frac{1}{2} & \frac{1}{2\sqrt{3}} & 0\\
    \frac{1}{2\sqrt{3}} & \frac{1}{6} &  0\\
    \frac{-\sqrt{3}}{\sqrt{2}} & \frac{-1}{\sqrt{2}} & \frac{1}{2}
    \end{array}\right].$$  Now $\rho(E_{\Phi,\Psi_1,\{1,2\}})=\rho(E_{\Phi,\Psi_1,\{1,3\}})=\rho(E_{\Phi,\Psi_1,\{2,3\}})=\max\left\{1,\frac{1}{3},0\right\}=1$, $\rho(E_{\Phi,\Psi_1,\{4,5\}})=\max\left\{1,0\right\}=1$ and $\rho(E_{\Phi,\Psi_1,\{1,4\}})=\rho(E_{\Phi,\Psi_1,\{1,5\}})=\rho(E_{\Phi,\Psi_1,\{2,4\}})=\rho(E_{\Phi,\Psi_1,\{2,5\}})=\rho(E_{\Phi,\Psi_1,\{3,4\}})=\rho(E_{\Phi,\Psi_1,\{3,5\}})=\max\left\{\frac{2}{3},\frac{1}{3},0\right\}=\frac{2}{3}$. Thus, $\rho_{\Phi,\Psi_1}^{(2)}=\max\left\{1,\frac{2}{3}\right\}=1$. Hence, $\Psi_1$ is an $SOD$-frame of $\Phi$ for $2$-erasures. Therefore, the canonical dual is a non-unique $SOD$-frame of $\Phi$ for $2$-erasures.}
\end{example}

	\section*{Acknowledgments}
	Aniruddha Samanta expresses thanks to the National Board for Higher Mathematics (NBHM), Department of Atomic Energy, India, for providing financial support in the form of an NBHM Post-doctoral Fellowship (Sanction Order No. 0204/21/2023/R\&D-II/10038). The second author also acknowledges excellent working conditions in the Theoretical Statistics and Mathematics Unit, Indian Statistical Institute Kolkata. 
	
	\mbox{}

\end{document}